\documentclass[12pt,]{article}
\usepackage[normal]{subfigure}
\usepackage{graphicx}
\usepackage{amssymb}
\usepackage{amscd}
\usepackage{amsmath}
\usepackage{amsfonts}
\usepackage{amsbsy}
\usepackage{fancyhdr}

\usepackage{epsfig}

\numberwithin{equation}{section}

\numberwithin{theorem}{section}
 
 \numberwithin{lemma}{section}
\date{}

\newcommand{\la}{\langle}
\newcommand{\ra}{\rangle}

\begin{document}
\date{}
\title{Limit cycle bifurcations from a nilpotent \\ focus or center of planar systems\thanks{The project supported  by  the  National
Natural Science Foundation of China (10971139) and the Slovenian
Research Agency.} }

\author{ Maoan Han$^{a,b}$\thanks {Corresponding author. {\em E-mail address}: mahan@shnu.edu.cn.},
Valery G. Romanovski$^{b,c}$\thanks {{\em E-mail
address}: valery.romanovsky@uni-mb.si.}\\
\small \emph{$^a$ Department of
Mathematics, Shanghai Normal University,}\\
\small \emph{ Shanghai
200234, P. R. China }
 \\
\small \emph{ $^b$ Faculty of Natural Science and Mathematics,}\\
\small \emph{University of Maribor,  SI-2000 Maribor, Slovenia}\\
\small \emph{$^c$ Center for Applied Mathematics and Theoretical
Physics, } \\
\small \emph{ University of Maribor, SI-2000 Maribor,
Slovenia }}

\date{}
\maketitle

 {\small \bf Abstract:}  { In this paper, we study the
analytical property of the Poincar\'e return map and the generalized
focal values  of an analytical planar system with a nilpotent focus
or center. Then we use the focal values and the map to study the
number of limit cycles of this kind of systems with parameters, and
obtain some new results on the lower and upper bounds of the maximal
number of limit cycles near the nilpotent focus or center.}
\hspace{4mm}

{\small \bf Keywords:}  Nilpotent focus; nilpotent center; limit cycle;  bifurcation. \ \hspace{4mm}

\section{Introduction  and main result}

Consider an analytic system of the form
$$ \dot x=y+X(x,y), \ \ \  \dot y=Y(x,y), \eqno(1.1)$$
where $X,\ Y=O(|x,y|^2)$ for $(x,y)$ near the origin.
The   following criterion for the existence of a center or a focus
at the origin of  (1.1) has been established in
 \cite{A,L}.

{\bf Theorem 1.1.} \  {\it Let $(1.1)$ have an isolated singular
point at the origin. Let
$$Y(x,F(x))=ax^{2n-1}+O(x^{2n}), \  a\neq 0,$$
$$ \frac{\partial X}{\partial x}(x,F(x))+\frac{\partial Y}{\partial
y}(x,F(x))=b x^{n-1}+O(x^n),$$
where $y=F(x)$ is the solution of the equation
$y+X(x,y)=0$ satisfying $F(0)=0$. Then the origin of $(1.1)$ is a
center or a focus if and only if $a$ is negative  and
$ b^2+4an<0.$}

Lyapunov \cite{L} also introduced the generalized polar
coordinates
$$x=r \ Cs(\theta), \ \ y=r^n\ Sn(\theta)$$
and the return map to give a way to find focal values in
solving the center-focus problem for (1.1), where $(Cs(t),Sn(t))$ is
the solution of the initial problem
$$\dot x=y, \ \dot y=-x^{2n-1}, \ \ (x(0),y(0))=(1,0).$$
  Sadovski \cite{Sad1} (see also \cite{ALS}) and
 Moussu \cite{M} investigated the problem  using
 Lyapunov function (Lyapunov constants) and normal form,  respectively.
 Then different ways of obtaining the focal values, Lyapunov
constants or their equivalent values and the bifurcation method of
local limit cycles were further given by Chavarriga, Giacomini, Gine
\& Llibre \cite{C}, Alvarez \& Gasull \cite{AG1,AG2} and Liu \& Li
\cite{LL1,LL2,LL3,LL4}. From Takens \cite{T} we know that (1.1) can
be formally transformed into  a formal normal form
$$\dot x=y, \  \dot y=-g(x)-yf(x), \eqno(1.2)$$
where $g(x)=ax^m+O(x^{m+1})$, $m\geq 2$  (the system (1.2) is a generalized Li\'enard system).
 Then,  Str\'o\.zyna \& \.Zo\l\c{a}dek
\cite{SZ} proved that this formal normal form can be achieved
through an analytic change of variables. Thus, if (1.1) has a center
or focus at the origin, then it can be changed into  (1.2)
with
$$g(x)=x^{2n-1}(a_{2n-1}+O(x)), \ n\geq 2, \ a_{2n-1}>0.
\eqno(1.3)$$ From Alvarez \& Gasull \cite{AG2} we see that under
(1.3) through a change of variables $x$ and $t$  of the form
$$u=[2n\int_0^xg(x)dx]^{\frac{1}{2n}}({\rm sgn }x)\equiv u(x),\
\frac{dt}{dt_1}=\frac{u^{2n-1}(x)}{g(x)}$$ the  system (1.2) can be
changed into
$$\dot x=y, \  \dot y=-x^{2n-1}-y\bar f(x), \eqno(1.4)$$
where
$$\begin{array}{c}
\bar f(x)=x^{2n-1}f(u^{-1}(x))/g(u^{-1}(x)), \ n\geq 2,\\
u(x)=[2n\int_0^xg(x)dx]^{\frac{1}{2n}}({\rm sgn}
x)=(a_{2n-1})^{\frac{1}{2n}}(x+O(x^2)).\end{array} \eqno(1.5)$$ Then,
by Theorem 1.1  system (1.4) has a center or a focus at the
origin if and only if the function $\bar f$ given in (1.5) satisfies
$$\bar f(x)=\sum_{j\geq n-1}b_jx^j,
\ \ b_{n-1}^2-4n<0. \eqno(1.6)$$ By Filippov's  theorem (see e. g. Ye  et
al. \cite{Y}) under (1.6) the  system (1.4) has a stable
(unstable) focus at the origin if there exists an integer $l$ with
$2l\geq n-1$ such that
$$ b_{2l}>0(<0), \  b_{2j}=0 \ {\rm for} \ j<l, \eqno(1.7)$$
and it has a center at the origin if $b_{2j}=0$ for all $2j\geq
n-1$.

Passing to  the generalized polar coordinate $(x,y)=(r  Cs(\theta), r^n
Sn(\theta))$ we obtain from  the system (1.4) the equation
$$ \frac{dr}{d\theta}=\frac{\sum_{j\geq
n-1}b_j(Sn(\theta))^2(Cs(\theta))^jr^{2-n+j}}{1+\sum_{j\geq
n-1}b_jSn(\theta)(Cs(\theta))^{j+1}r^{1-n+j}}. \eqno(1.8)$$
The function on the right hand side of (1.8) is
 periodic of the period
$T=2\sqrt{\frac{\pi}{n}}\Gamma(\frac{1}{2n})/\Gamma(\frac{n+1}{2n}).$
Let $r(\theta,r_0)$ denote the solution of (1.8) with the initial
value $r(0)=r_0$. Then
$$r(T,r_0)=\sum_{j\geq 1}V_jr_0^j.$$
 Alvarez \& Gasull \cite{AG2} called the constant $V_k$
the $k$th generalized Lyapunov constant of (1.8) assuming
$V_1=1,V_2=\cdots=V_{k-1}=0$. They  also
studied the normal form (1.4) and proved the following theorem.

{\bf Theorem 1.2.} \ {\it Let $(1.6)$ and $(1.7)$ be satisfied. Then

$(1)$ \ $V_1=\exp(\frac{-2b_{n-1}\pi}{n\sqrt{4n-b_{n-1}^2}})$   if
$2l=n-1$;

$(2)$ \ $V_1=1$, $V_j=0$ for $1<j<2-n+2l$, and
$V_{2-n+2l}=-K_lb_{2l}$ if either $b_{n-1}=0$ or $b_{n-1}\neq 0$ and
$n$ is even, where $K_l$ is a positive constant.}

  For the case of $n=2$ Liu and Li \cite{LL1} introduced a different
  generalized polar coordinates of the form $x=r\cos\theta$,
  $y=r^2\sin\theta$ to change  (1.1) into the form
  $$\frac{dr}{dt}=R(\theta,r),\ \frac{d\theta}{dt}=Q(\theta, r),$$
  assuming the origin is a center or a focus. Let $\tilde r(\theta,h)$
  denote the solution of the $2\pi$-periodic system
  $$\frac{dr}{d\theta}=\frac{R(\theta,r)}{Q(\theta,r)}$$
  satisfying $\tilde r(0)=h$. Note that the initial value problem is
well-defined also for negative $h$.

For analytic functions $\phi, \phi_1,\dots,\phi_k  $ defined on a domain $D$ we will write $\phi=O(| \phi_1,\dots,\phi_k |) $ if there are analytic functions $\psi_1,\dots, \psi_k$ on
$D$, such that  $\phi=\psi_1 \phi_1+\dots+\psi_k \phi_k   $ on $D$. Liu and Li \cite{LL1} found the
  following facts.

  {\bf Theorem 1.3.}  \ {\it  Consider the system  $(1.1)$. Let the conditions of
  Theorem $1.1$ be satisfied with $n=2$ $($or $m=3$ $)$ such that the origin
  is a center or a focus. Then,

  $(1)$ \ $\tilde r(\theta,-\tilde r(\pi,h))=-\tilde
  r(\pi-\theta,h)$;

  $(2)$ \ $\Delta (h)=\tilde r(-2\pi,h)-h=\sum_{k\geq 2} v_kh^k$,
  where
  $$v_{2k+1}=O(| v_2,v_4,\cdots,v_{2k}|), \ k\geq 1;$$

  $(3)$ \ the origin is a stable  $($unstable$)$ focus if
  $$v_{2k}<0(>0),\ {\it and } \ v_{2j}=0 \ {\it  for } \ j<k.$$
  In the latter  case  the origin is called   a $k$th order weak focus
  of $(1.1)$.}

 Liu and Li \cite{LL1} also gave some new methods to compute the
 focus values $v_2,v_4,\cdots,v_{2k}$, or  equivalent values,
 and studied the problem of limit cycle bifurcations near the
 origin, finding a new phenomenon: a node can generate a limit cycle when its stability changes.

 In this paper we study the problem of limit cycle bifurcations near
 the origin for the  analytic system
 $$\dot x=y+X(x,y,\delta), \ \dot y=Y(x,y,\delta),\eqno(1.9)$$
where $\delta =(\delta_1,\dots, \delta_m)\in D \subset \mathbb{R}^m$ with $D$ compact, and $X,\
Y=O(|x,y|^2)$ for $ |x|$ small and $\delta \in D$. Let
$y=F(x,\delta)$ be the solution of the equation $y+X(x,y,\delta)=0$.
We define the following two functions:
$$g(x,\delta)=-Y(x,F(x,\delta),\delta),\ f(x,\delta)=-\left[\frac{\partial
X}{\partial x}(x,F(x,\delta),\delta)+\frac{\partial Y}{\partial
y}(x,F(x,\delta),\delta)\right].\eqno(1.10)$$ By Theorem 1.1,
if
$$ g(x,\delta)=x^{2n-1}(a_{2n-1}(\delta)+O(x)), \ n\geq 2,\
a_{2n-1}(\delta)>0,\eqno(1.11)$$ $$ f(x,\delta)=\sum_{j\geq
n-1}b_j(\delta)x^j, \ \
b_{n-1}^2(\delta)-4na_{2n-1}(\delta)<0,\eqno(1.12)
$$ then the origin is a  center or a  focus of (1.9) for all $\delta \in
D$.

Let us define a Poincar\'e return map for the plane system (1.9).
For each $\delta \in D$ and $x_0\neq 0$ with $|x_0|$ small consider
the solution $(x(t,x_0,\delta),y(t,x_0,\delta))$ of (1.9) with the
initial condition $(x(0),y(0))=(x_0,F(x_0,\delta))$. Then there is a
unique least positive number $\tau=\tau(x_0,\delta)>0$ such that
$y(\tau,x_0,\delta)=F(x(\tau,x_0,\delta),\delta)$ and
$x_0x(\tau,x_0,\delta)>0$. See Figure 1 for $x_0>0$ small.

\begin{center}
\includegraphics[bb=0 145 592 730,scale=0.33]{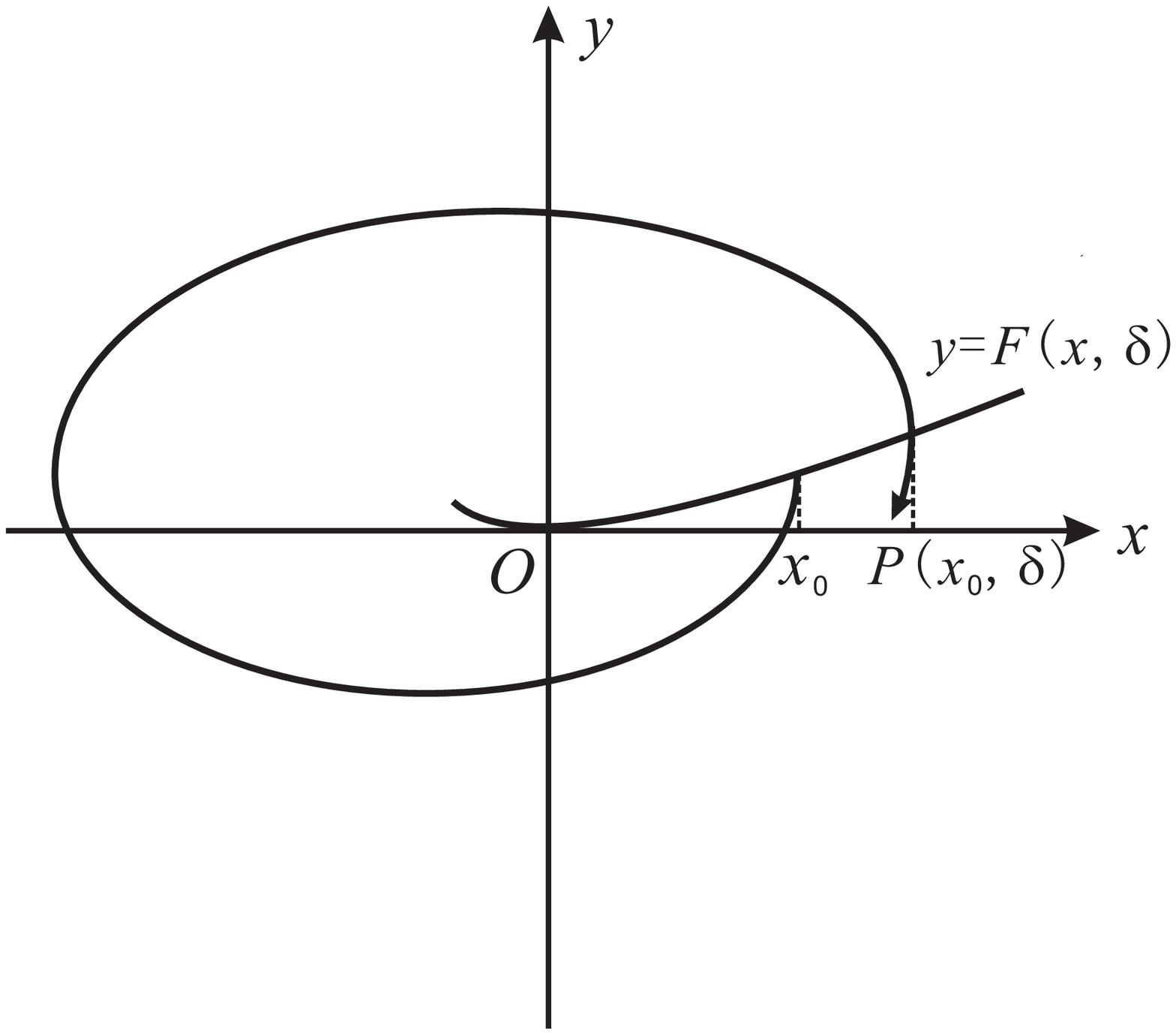}\\
{Figure 1. The Poincar\'{e} map of (1.9) with $x_0>0$.}
\end{center}

 Thus, the Poincar\'e return map is
defined as
$$P(x_0,\delta)=\left\{\begin{array}{ll} x(\tau,x_0,\delta), & 0<|x_0|<\varepsilon_0,\\
0,& x_0=0\end{array}\right.$$ where $\varepsilon_0$ is a small
positive constant. Evidently, the function is uniquely defined, and
it is continuous at $x_0=0$ under (1.11) and (1.12). Moreover, (1.9)
has a periodic orbit near the origin if and only if the map has two
fixed points near zero: one positive and the other one negative. For
the analytical property of this function at $x_0=0$, we have the
following theorem.

{\bf Theorem 1.4.} \  {\it Let $(1.9)$ satisfy $(1.11)$ and $(1.12)$
for all $\delta \in D$. Then there is a unique  analytic function
$\bar P(x_0,\delta)$ in $x_0$ at $x_0=0$, satisfying $\frac{\partial
\bar P}{\partial x_0}(0,\delta)>0$ and having the expansion
$$\bar d(x_0,\delta)=\bar P(x_0,\delta)-x_0=\sum_{j\geq 1}v_j(\delta)x_0^j\eqno(1.13)$$ for
$|x_0|$ small, such that

$(1)$ \ if $n$ is odd, then  $P(x_0,\delta)=\bar
P(x_0,\delta)$ for all $|x_0|$ small;

$(2)$ \ if $n$ is even, then for all $|x_0|$ small
$$ P(x_0,\delta)=\left\{\begin{array}{ll} \bar P(x_0,\delta) \ \  {\rm
for } \ x_0>0,\\ \bar P^{-1}(x_0,\delta) \ {\rm for } \
x_0<0,\end{array} \right. $$ where $\bar P^{-1}$ denotes the inverse
of $\bar P$ in $x_0$.

Hence, the system $(1.9)$ has a periodic orbit near the origin if
and only if the analytic function $\bar d$ defined in $(1.13)$ has
two zeros in $x_0$ near $x_0=0$, among which   one is  positive and the other one is
negative. The function $\bar d$ is called the succession function or the
bifurcation function of $(1.9)$.}

The above theorem tells us that the function $ P(x_0,\delta)$ is
analytic in $x_0$ at $x_0=0$ as $n$ is odd, and not analytic in
$x_0$ at $x_0=0$ as $n$ is even unless the origin is a center (in
this case, $P$ is the identity).

For the property of the coefficients $v_j$ in (1.13) we have further

{\bf Theorem 1.5.} \ {\it Let $(1.9)$ satisfy $(1.11)$ and $(1.12)$
for all $\delta \in D$. Then

$(1)$ \ For $n$ odd we have $
 v_{2k}=O(|v_1,v_3,\cdots,v_{2k-1}|), \ k\geq 1.$

 $(2)$ \ For $n$ even we have  $v_1=0$, $
 v_{2k+1}=O(|v_2,v_4,\cdots,v_{2k}|), \ k\geq 1.$
}

Define $p_n=[1+(-1)^n]/2$. Then the conclusions of the above theorem
can be written  uniformly as
$$v_{2k+p_n}= O(|  v_{1+p_n},v_{3+p_n},\cdots,v_{2k-1+p_n}|), \ k\geq
1.$$

From the proof of the above theorem we  see that $v_{2k+p_n}$
depends on $v_{1+p_n}$, $v_{3+p_n}$, $\cdots$, $v_{2k-1+p_n}$
smoothly.
 Using the  theorem we derive the
following two  statements  on limit cycle bifurcations near the origin.

{\bf Theorem 1.6} (Bifurcation from Focus). \ {\it Let $(1.9)$
satisfy $(1.11)$ and $(1.12)$ for all $\delta \in D$. Denote
$p_n=[1+(-1)^n]/2$.

 $(1)$\ If there is an integer $k\geq 1$ such
that
$$\sum_{j=1}^{k+1}|v_{2j-1+p_n}(\delta)|>0  \ {\rm for  \ all} \ \delta
\in D,$$ then there exists a neighborhood $U$ of the origin such
that $(1.9)$ has at most $k$ limit cycles in $U$ for all $\delta \in
D$.

$(2)$\ If there is $\delta_0 \in D$ such that
$v_{2k+1+p_n}(\delta_0)\neq 0, $  and
$$\begin{array}{c} v_{2j-1+p_n}(\delta_0)= 0,
j=1,\cdots, k,\\
{\rm
rank}\frac{\partial(v_{1+p_n},v_{3+p_n},\cdots,v_{2k-1+p_n})}
{\partial(\delta_1,\delta_2,\cdots,\delta_m)}(\delta_0)=k,\end{array}
\eqno(1.14)$$ then for an arbitrary   sufficiently small neighborhood of the origin there
are some $\delta\in D$ near $\delta_0$ such that $(1.9)$ has exactly
$k$ limit cycles in the neighborhood.}

{\bf Theorem 1.7} (Bifurcation from Center). \ {\it Let $(1.9)$
satisfy $(1.11)$ and $(1.12)$ for all $\delta \in D$.
 Assume that there exist $\delta_0
\in D$ and an integer $k\geq 1$ such that $(1.14)$ is satisfied. If
the origin is a center of $(1.9)$ as $v_{2j-1+p_n}(\delta)= 0,
j=1,\cdots, k,$ then there exists a neighborhood $U$ of the origin
such that $(1.9)$ has at most $k-1$ limit cycles in $U$ for all
$\delta\in D$ near $\delta_0$, and also, for an arbitrary
 sufficiently small
neighborhood of the origin there are some $\delta\in D$ near
$\delta_0$ such that $(1.9)$ has exactly $k-1$ limit cycles in the
neighborhood.}

The theorem means that the cyclicity of the system at the point $\delta_0$ is equal to
$k-1$.

%More precisely bifurcations from a center can be described by an analog of Christopher's theorem %\cite{Ch}.

Now, different from \cite{AG2} and \cite{LL1}--\cite{LL4}, we give
the following new and more reasonable definition.

{\bf Definition 1.1.} \ We call $v_{2k+1+p_n}(\delta)$ the
generalized focal values of order $k$   of
(1.9) at the origin.

 By Theorem 1.6, we see that a nilpotent focus of order
 $k$ generates at most $k$ limit cycles under perturbations as long
 as the perturbations always satisfy $(1.11)$ and $(1.12)$.

The generalized focal values
$v_{1+p_n},v_{3+p_n},\cdots,v_{2k+1+p_n},\cdots$  can be calculated
 using the normal form  of system  (1.9). We will give a method
 how to do it. By Str\'o\.zyna and \.Zo\l\c{a}dek \cite{SZ} we know
that (1.9) has the following analytic normal form:
$$\dot x=y, \  \dot y=-g(x,\delta)-yf(x,\delta). \eqno(1.15)$$
We remark that here $f$ and $g$ in (1.15) may be different from ones
given by (1.10). As before, let $\delta \in D \subset \mathbb{R}^m$
with $D$ compact. Also, suppose for $ |x|$ small the function
$g(x,\delta)$ satisfies (1.11).
 Define
 $$F(x,\delta)=\int_0^xf(x,\delta)dx,\ \
 G(x,\delta)=\int_0^xg(x,\delta)dx.$$ It is easy to see that the
 equation $G(x,\delta)=G(y,\delta)$ for $xy<0$ defines a unique
 analytic function $y=\alpha(x,\delta)=-x+O(x^2)$. Introduce
 $$F(\alpha(x,\delta),\delta)-F(x,\delta)=\sum_{j\geq 1}B_j(\delta)x^j.\eqno(1.16)$$

By Theorem 1.1, if (1.15) satisfies (1.11) and (1.12) then it has a
center or focus at the origin. Thus, under (1.11) and (1.12) the
Poincar\'e return map for (1.15) is well defined near the origin.

{\bf Theorem 1.8.} \  {\it  Let $(1.15)$ satisfy $(1.11)$ and
$(1.12)$ for  all  $\delta \in D$. Then for $x_0>0$ small,  the
Poincar\'e return map $P(x_0,\delta)$ has the form
$$P(x_0,\delta)-x_0=\sum_{j\geq
0}v_{2j+1+p_n}(\delta)x_0^{2j+1+p_n}(1+P_j^*(x_0,\delta)),$$ where
$P_j^*(x_0,\delta)=O(x_0)$,
$$\begin{array}{l}
v_{1+p_n}(\delta)=K_l^*B_{2l+1}(\delta)+(1-p_n)O( B_{2l+1}^2),\\
v_{2j+1+p_n}(\delta)=K_{l+j}^*
B_{2l+2j+1}(\delta)+O(|B_{2l+1},B_{2l+3},\cdots,B_{2l+2j-1}|), \
j\geq 1,\end{array} \eqno(1.17)$$ $l=[n/2]$, and $K_{l+j}^*,\ j\geq
0$ are positive constants. Thus, Theorems $1.6$ and $1.7$ hold if
$v_{2j+1+p_n}$ is replaced by $B_{2l+2j+1}$, $j\geq 0$. }

Let $$ f(x,\delta)=\sum_{j\geq0}b_j(\delta)x^j. \eqno(1.18)$$

Then, we further have for (1.15)

{\bf Theorem 1.9.} \  {\it  Let $(1.15)$ satisfy $(1.11)$, $(1.16)$
and $(1.18)$ for all $\delta \in D$. Assume  there exist $\delta_0
\in D$ and $k\geq [n/2]$ such that
$$B_{2k+1}(\delta_0)<0(>0), \ B_{2j-1}(\delta_0)= 0, j=1,\cdots,
k.\eqno(1.19)$$ Let one of the following conditions be satisfied:

$(a)$ \ $n=2,$ and $$
 b_0(\delta_0)=0, \  b_{1}^2(\delta_0)-8a_{3}(\delta_0)<0;\eqno(1.20)$$

$(b)$    $n>2,$ $g(-x,\delta)=-g(x,\delta)$,
$f(-x,\delta)=f(x,\delta)$, and $$
 b_j(\delta_0)=0 \ {\it for} \ j=0,\cdots,n-2\ {\it  and}   \   b_{n-1}^2(\delta_0)-4na_{2n-1}(\delta_0)<0.\eqno(1.21)$$
 Then we have

$(1)$ \ For $\delta=\delta_0$ $(1.15)$ has a stable $($unstable$)$
focus at the origin.

$(2)$ \ If further $${\rm
rank}\frac{\partial(B_{1},B_{3},\cdots,B_{2k-1})}
{\partial(\delta_1,\delta_2,\cdots,\delta_m)}(\delta_0)=k,$$ then
for an arbitrary  sufficiently small  neighborhood of the origin there are some
$\delta\in D$ near $\delta_0$ such that $(1.15)$ has at least $k$
limit cycles in the neighborhood.}

From Theorems 1.4--1.8, it seems that under (1.11) and (1.12) we have
solved the problem of limit cycle bifurcation for generic systems.
Theoretically it is,
but in practice it is not. The reason is that in general we do not
know what is the transformation from (1.9) to its normal form
(1.15). Here we give a method to solve the problem  completely both
theoretically and in practice. It includes three main steps below.

First, under (1.11) and (1.12) by the  normal form theory (see, for instance, \cite{T}),
 for any integer $m>2n-1$ there is a change of
variables of the form
$$\left( \begin{array}{c} x\\ y\end{array}\right)=\left( \begin{array}{c} u\\ v\end{array}\right)+H_m(u,v,\delta),$$
where $H_m(u,v,\delta)=O(|u,v|)$ is a polynomial in $u,v$ of degree
at most $m$, such that it transforms (1.9) into (1.22)
(called the normal form of order $m$ of (1.9), or the Takens normal
form; we still use $(x,y)$ for the new variables $u,v$)
$$\dot x=y+X_{m+1}(x,y,\delta), \  \dot y=-g_m(x,\delta)-yf_{m-1}(x,\delta)+Y_{m+1}(x,y,\delta), \eqno(1.22)$$
where
$$g_m(x,\delta)=\sum_{j=2n-1}^ma_j(\delta)x^j, \ \ f_{m-1}(x,\delta)=\sum_{j=n-1}^{m-1}b_j(\delta)x^j
$$ with $a_{2n-1}(\delta)>0$ and $
b_{n-1}^2(\delta)-4na_{2n-1}(\delta)<0$, and $X_{m+1}(x,y,\delta)$,
$Y_{m+1}(x,y,\delta)$ being  analytic functions satisfying $X_{m+1},
\  Y_{m+1}=O(|x,y|^{m+1})$. Here, we should mention that the
functions $g_m$ and $f_{m-1}$ depend only on the terms of degree at
most $m$ of the expansions of the functions $X$ and $Y$ in (1.9) at
the origin.

The Poincar\'e maps of (1.9) and (1.22) are essentially the  same. We can
suppose that the Poincar\'e map of (1.22) is $P(x_0,\delta)$ having
the expansion
$$ P(x_0,\delta)-x_0=\sum_{j\geq 1}v_j(\delta)x_0^j \eqno(1.23)$$
for $x_0>0$ small.

Second, truncating the higher order terms in (1.22) we obtain the
following polynomial system of degree $m$
$$\dot x=y, \  \dot y=-g_m(x,\delta)-yf_{m-1}(x,\delta). \eqno(1.24)$$
 In practice, for given system (1.9) it is not difficult to find the
 corresponding system (1.24). For (1.24) we can further use Theorem
 1.8 to find its focal values at the origin up to any large order. Let
 $P_m(x_0,\delta)$ denote the Poincar\'e map of (1.24). It has the
 expansion
$$ P_m(x_0,\delta)-x_0=\sum_{j\geq 1}\bar v_j(\delta)x_0^j \eqno(1.25)$$
for $x_0>0$ small.

Third, we want to use $\bar v_j(\delta)$ for  $ v_j(\delta)$. Here,
a problem we would like to solve is the following: For any given
$k>1$ find $m>2n-1$ such that $v_j(\delta)=\bar v_j(\delta)$ for $
1\leq j \leq k.$ The following theorem gives an answer.

{\bf Theorem 1.10.} \  {\it   Consider $(1.22)$ and $(1.24)$. Then
for any integer $k\geq 1$, if $m\geq (k+2)n-2$ then
$$v_j(\delta)=\bar v_j(\delta) \ \   {\rm for} \  1\leq j \leq k\, n. \eqno(1.26)$$}

Therefore, we have

{\bf Corollary 1.1.} \ {\it Under $(1.11)$ and $(1.12)$ for any
integer $k\geq 1$  for $(1.9)$ the coefficients
$v_1,v_2,\cdots,v_{kn}$ in $(1.13)$ depend only on the terms of
degree at most $(k+2)n-2$ of the expansions of the functions $X$ and
$Y$ at the origin.}

Obviously, in the case of $n=1$ (the elementary case), the above
conclusion is a well-known results.

 We organize the paper as follows. In
section 2 we first give preliminary lemmas. In section 3 we prove
our main results. In section 4 we provide some application examples.

\section{ Preliminaries }

 Consider (1.9). In this section we will always suppose that  (1.11) and (1.12) are satisfied.
 Introducing a new variable $v=y-F(x,\delta)$ we can
 obtain from (1.9)  (reusing $y$ for $v$)
 $$\begin{array}{l}
 \dot x=y(1+Z_1(x,y,\delta)),\\
 \dot y=-g(x,\delta)-yf(x,\delta)+y^2Z_2(x,y,\delta),\end{array}
 \eqno(2.1)$$
 where the functions $f$ and $g$ are given by (1.10), and $Z_1$ and
 $Z_2$ are analytic functions near the origin with
 $Z_1(x,y,\delta)=O(|x,y|)$. In the discussion below we will often
 omit $\delta$ for convenience.
 As in Liu and Li \cite{LL4} we will make a change of variables to (2.1)
 using the generalized polar coordinates
 $$x=r\cos\theta,\ \ y=r^n\sin\theta, \ r>0.\eqno(2.2)$$

 {\bf Lemma 2.1.} \ {\it Let $(1.11)$ and $(1.12)$ be satisfied. Then the
 transformation $(2.2)$ carries $(2.1)$ into the form
 $$\begin{array}{l}
 \dot \theta=S(\theta,r)=\frac{r^{n-1}}{H(\theta)}[S_0(\theta)+O(r)],\\
  \dot
  r=R(\theta,r)=\frac{r^{n}}{H(\theta)}[R_0(\theta)+O(r)],\end{array}\eqno(2.3)
  $$
  where $S$ and $R$ are $2\pi$-periodic in $\theta$, and satisfy
  $$S(\pi+(-1)^{n-1} \theta,-r)=(-1)^{n-1}S(\theta,r),\ R(\pi+(-1)^{n-1}
  \theta,-r)=-R(\theta,r),\eqno(2.4)$$ and $H(\theta)=\cos^2\theta+n\sin^2\theta>0$,
  $$S_0(\theta)=-[n\sin^2\theta+b_{n-1}\cos^n\theta\sin\theta+a_{2n-1}\cos^{2n}\theta]<0,$$
  $$R_0(\theta)=\cos\theta\sin\theta(1-a_{2n-1}\cos^{2n-2}\theta-b_{n-1}\sin\theta\cos^{n-2}\theta).$$
}

{\bf Proof.}\ From (2.2) we have
$$\dot x=\cos\theta\dot r-r\sin\theta\dot \theta,\ \dot
y=nr^{n-1}\sin\theta\dot r+r^n\cos\theta\dot \theta.$$ We
solve the above equations for $\dot \theta$ and $\dot r$,  and obtain (2.3) with
$$S(\theta,r)=\frac{\cos\theta\dot y-nr^{n-1}\sin\theta \dot
x}{r^n(\cos^2\theta+n\sin^2\theta)},$$
$$R(\theta,r)=\frac{\sin\theta\dot y+ r^{n-1}\cos\theta\dot
x}{r^{n-1}(\cos^2\theta+n\sin^2\theta)}.$$ Then noting that
$$\cos(\pi\pm \theta)=-\cos\theta,\ \sin(\pi\pm
\theta)=\mp\sin\theta $$ and that (2.2) is invariant as $(\theta,r)$
is replaced by $(\pi+(-1)^{n-1} \theta,-r)$ one can prove (2.4)
easily. The other conclusions are direct. This ends the proof.

By (2.3) and (2.4) we obtain the following analytic $2\pi$-periodic
equation
$$\frac{dr}{d\theta}=\bar R(\theta,r),\eqno(2.5)$$
where
$$\begin{array}{c}
\begin{array}{cl}
\bar R(\theta,r)&=\displaystyle r\frac{\sin\theta\dot y+
r^{n-1}\cos\theta\dot
x}{\cos\theta\dot y-nr^{n-1}\sin\theta \dot x}\\
&=r[R_0(\theta)/S_0(\theta)+O(r)],\end{array}\\
 \bar
R(\pi+(-1)^{n-1}
  \theta,-r)=(-1)^nR(\theta,r).\end{array}\eqno(2.6)$$

Let $r(\theta,h)$ denote the solution of (2.5) with the initial
value $r(0)=h$. For properties of the solution we have

{\bf Lemma 2.2.}\ {\it The solution $r(\theta,h)=O(h)$ is analytic
in $(\theta,h)$ for $|h|$ small, and satisfies

$(1)$\ $r(\theta,-r(\pi,h))=-r(\pi+(-1)^{n-1}\theta,h)$;

$(2)$\ $r(\theta\pm 2\pi,h)=r(\theta,r(\pm 2\pi,h))$.}

{\bf Proof.} \ Let $\tilde r(\theta)=-r(\pi+(-1)^{n-1}\theta,h)$.
Then by (2.5) and (2.6) we have
$$\begin{array}{rl}
\frac{d\tilde r}{d\theta}&=(-1)^n\bar
R(\pi+(-1)^{n-1}\theta,r(\pi+(-1)^{n-1}\theta,h))\\
&=(-1)^n\bar
R(\pi+(-1)^{n-1}\theta,-\tilde r(\theta))\\
&=\bar R(\theta,\tilde r(\theta)).\end{array}$$ This means that
$\tilde r(\theta)$ is also a solution of (2.5). Then the first
conclusion follows by the uniqueness of initial problem. The second
one follows in the same way. This completes the proof.

Further we have

{\bf Lemma 2.3.} \ {\it Let $P(x_0,\delta)$ be the Poincar\'e return
map of $(1.9)$ defined in section $1$. Then  for $|x_0|>0$ small we
have $P(x_0,\delta)=r(-2\pi,x_0)$ for $x_0>0$, and
$P(x_0,\delta)=r((-1)^n2\pi,x_0)$ for $x_0<0$.}

{\bf Proof.} \ First, it is easy to see that (1.9) and (2.1) have
the same Poincar\'e return map $P(x_0,\delta)$. Then, noting that
 $\dot
{\theta}<0$ for $r>0$ small by (2.3), by the definition of $P$ and
(2.2) we can see that
$$P(x_0,\delta)=x(\tau,x_0)=r(-2\pi, x_0)$$
for $x_0>0$ small. Now consider the case of  $x_0<0$. Let
$r^*(\theta,h)$ denote the solution of (2.5) satisfying
$r^*(\pi)=h$. Then we have similarly
$$P(x_0,\delta)=x(\tau,x_0)=-r^*(-\pi, -x_0),$$
since under (2.2) the points $(x_0,0)$ and $(P(x_0,\delta),0)$ on
the $(x,y)$-plane correspond to the points $(\pi,-x_0)$ and
$(-\pi,-P(x_0,\delta))$ on the $(\theta,r)$-plane respectively.

Further, by Lemma 2.2(1) we have
$$r^*(\theta,-h)=-r(\pi-\theta,h) \ {\rm for} \ n \ {\rm even},\eqno(2.7)$$
and
$$r^*(\theta,-r(2\pi,h))=-r(\pi+\theta,h) \ {\rm for} \ n \ {\rm odd}. \eqno(2.8)$$
Noting that by Lemma 2.2(2), $x_0=r(2\pi,h)$ if and only if
$h=r(-2\pi,x_0)$, we see that  (2.8) becomes
$$r^*(\theta,-x_0)=-r(\pi+\theta,r(-2\pi,x_0)) \ {\rm for} \ n \ {\rm odd}. \eqno(2.9)$$
Therefore, for $x_0<0$ by (2.7) and (2.9)
$$P(x_0,\delta)=\left\{\begin{array}{l}r(2\pi,x_0) \ \ {\rm for}\ n\
 {\rm even},\\
r(-2\pi,x_0) \ {\rm for}\ n \ {\rm odd}. \end{array}\right.$$ This
ends the proof.

{\bf Lemma 2.4.} \ {\it Let $d(x_0,\delta)=P(x_0,\delta)-x_0$. Then
there exists an analytic function $K(h,\delta)$ for $|h|$ small with
$K(0,\delta)=\frac{\partial r}{\partial x_0}(\pi,0)>0$ such that
$$d(\tilde x_0,\delta)=-K(x_0,\delta)d(x_0,\delta)\eqno(2.10)$$
for $x_0>0$ small, where
 $\tilde x_0=-r(\pi,x_0)$. }

{\bf Proof.} \  By Lemma 2.2, we have
$$r((-1)^n2\pi,\tilde x_0)=-r(-\pi,x_0)=-r(\pi,r(-2\pi,x_0)).$$
Hence, by Lemma 2.3 for $x_0>0$
$$\begin{array}{cl}
d(\tilde x_0,\delta)&=r((-1)^n2\pi,\tilde x_0)-\tilde x_0\\
&=-r(\pi,r(-2\pi,x_0))+r(\pi,x_0)\\
&=-K(x_0,\delta)[r(-2\pi,x_0)-x_0]\\
&=-K(x_0,\delta)d( x_0,\delta),\end{array}$$ where
$$K(x_0,\delta)=\int_0^1\frac{\partial r}{\partial
x_0}(\pi,x_0+s(r(-2\pi,x_0)-x_0))ds.$$ It is obvious that $K$ is
analytic for $|x_0|$ small and $K(0,\delta)=\frac{\partial
r}{\partial x_0}(\pi,0)>0$. This completes the proof.

\section{Proof of the main results}

In this section we prove our main results presented in  Theorems 1.4--1.10.

{\bf Proof of Theorem 1.4.} \ We take $\bar
P(x_0,\delta)=r(-2\pi,x_0)$ for $|x_0|$ small. Then by Lemma 2.2
$\bar P$ is analytic. Note that by Lemma 2.2, $r(2\pi,x_0)$ is the
inverse of $r(-2\pi,x_0)$ in $x_0$. Then Theorem 1.4 follows
directly from Lemma 2.3. The proof is complete.

{\bf Proof of Theorem 1.5.} \  There are two cases to consider
separately.

{\bf Case A:}\ $n$ odd. By (1.13) and Theorem 1.4(1), we have
$$d(x_0,\delta)=\bar d(x_0,\delta)=\sum_{j\geq
1}v_j(\delta)x_0^j\eqno(3.1)$$ for all $|x_0|$ small.

By Lemma 2.4, we can suppose
$$K(x_0,\delta)=\sum_{j\geq 0}k_jx_0^j,\ \ \tilde x_0=-r(\pi,x_0)=\sum_{j\geq
1}l_jx_0^j,\eqno(3.2)$$ where $k_0>0$, $l_1=-k_0$. Substituting
(3.1) and (3.2) into (2.10), we obtain
$$\sum_{j\geq
1}v_j(\sum_{i\geq1}l_ix_0^i)^j=-\sum_{
i\geq0,j\geq1}k_iv_jx_0^{i+j}.$$ Comparing the coefficients of the
terms $x_0^{2}$, $x_0^{4}$ and $x_0^{2j}$ on both sides yields
$$\begin{array}{l}
v_2l_1^2+v_1l_2=-(k_0v_2+k_1v_1),\\
v_4l_1^4+3v_3l_1^2l_2+v_2(l_2^2+2l_1l_3)+v_1l_4=-(k_0v_4+k_1v_3+k_2v_2+k_3v_1),\\
\cdots \cdots \cdots \cdots \\
v_{2j}l_1^{2j}+v_{2j-1}L_{1,j}(l_1,l_2)+\cdots+v_2L_{2j-2,j}(l_1,l_2,\cdots,l_{2j-1})+v_1l_{2j}=-\sum_{i=0}^{2j-1}k_iv_{2j-i},\\
\cdots \cdots \cdots \cdots \end{array}$$ where
$L_{i,j}(l_1,l_2,\cdots,l_{i+1})$, $i=1,2,\cdots,2j-2$, are all
polynomials. Then from the above equations we obtain
$$
 v_{2k}=O(| v_1,v_3,\cdots,v_{2k-1}|), \ k\geq 1.$$

{\bf Case B:}\ $n$ even. By (1.13) and $x_0=\bar P^{-1}(\bar
P(x_0,\delta),\delta)$ we can find
$$\bar P^{-1}(x_0,\delta)=\tilde v_1x_0+\tilde v_2x_0^2+\tilde
v_3x_0^3+\cdots, \eqno(3.3)$$ where
$$\begin{array}{l}
\tilde v_1=(v_1+1)^{-1},\\
\tilde v_2=-v_2(v_1+1)^{-3},\\
\cdots \cdots \\
\tilde
v_j=-v_j(v_1+1)^{-(j+1)}+L_j(v_2,v_3,\cdots,v_{j-1}),\\
\cdots \cdots \end{array} \eqno(3.4)$$ where each $L_j$ is a
polynomial of degree at least 2. Now we suppose $x_0>0$. Then (3.1)
holds by Theorem 1.4. Further, noting that  $\tilde x_0<0$ by Theorem 1.4
again
$$d(\tilde x_0,\delta)=P(\tilde x_0,\delta)-\tilde x_0=\bar P^{-1}(\tilde x_0,\delta)-\tilde
x_0.\eqno(3.5)$$ Then,  inserting (3.1), (3.2), (3.3) and (3.5) into
(2.10) we obtain
$$\begin{array}{l}
(\tilde v_1-1)l_1=-k_0v_1,\\
(\tilde v_1-1)l_2+\tilde v_2l_1^2=-(k_0v_2+k_1v_1),\\
\cdots \cdots \\
(\tilde v_1-1)l_j+L_j(\tilde v_2,\tilde v_3,\cdots,\tilde v_{j-1})+\tilde v_jl_1^j=-(k_0v_j+k_1v_{j-1}+\cdots+k_{j-1}v_1),\\
\cdots \cdots \end{array}\eqno(3.6)$$ where $$L_j(\tilde v_2,\tilde
v_3,\cdots,\tilde v_{j-1})\in \la \tilde v_2,\tilde v_3,\cdots,\tilde
v_{j-1}\ra.$$ Finally, noting that  $l_1=-k_0$ and  substituting (3.4) into
(3.6) we easily see that
 $$v_{2j+1}=O(| v_2, v_4,\cdots,
v_{2j}|), \ j\geq 1.$$ This ends the proof.

{\bf Proof of Theorem 1.6.} \ For the first part, suppose the
conclusion is not true. Then there exists a sequence $\{\delta_m\}$
in $D$ such that for $\delta=\delta_m$ (1.9) has $k+1$ limit cycles
$L_{m,1},L_{m,2},\cdots,L_{m,k+1}$ which approach the origin as
$m\rightarrow\infty$. Then by Theorem 1.5, the function $\bar
d(x_0,\delta_m)$ has $2k+2$ non-zero roots in $x_0$ which approach
zero as $m\rightarrow\infty$.

Since $D$ is compact, we can assume $\delta_m\rightarrow\delta_0\in
D$ as $m\rightarrow\infty$. By our assumption,
$\sum_{j=1}^{k+1}|v_{2j-1+p_n}(\delta_0)|>0$. Thus, for some $1\leq
l\leq k+1$,
$$v_{2l-1+p_n}(\delta_0)\neq0, \ v_{2j-1+p_n}(\delta_0)=0\ {\rm
for}\ 1\leq j\leq l-1.$$ Therefore, by (1.13) and Theorem 1.5, we
have
$$\bar d(x_0,\delta_0)=v_{2l-1+p_n}(\delta_0)x_0^{2l-1+p_n}+O(x_0^{2l+p_n}).$$
Note that $\bar d(0,\delta)=0$. It follows from Rolle's theorem that
for some $\varepsilon_0>0$ the function $\bar d(x_0,\delta)$ has at
most $2l-2+p_n$ non-zero roots in $(-\varepsilon_0,\varepsilon_0)$
for all $|\delta-\delta_0|<\varepsilon_0$. We have proved that the
function $\bar d(x_0,\delta_m)$ has $2k+2$ non-zero roots which
approach zero as $m\rightarrow\infty$. It then follows that
$2k+2\leq 2l-2+p_n$, contradicting to $2l-2+p_n\leq 2k+p_n\leq
2k+1$. The first conclusion follows.

For the second one, by Theorem 1.5, the function $\bar d$ can be
written as
$$\bar d(x_0,\delta)=\sum_{j\geq 1}v_{2j-1+p_n}(\delta)x_0^{2j-1+p_n}(1+P_j(x_0,\delta)),\eqno(3.7)$$
where $P_j(0,\delta)=0$. Like   in \cite{CJ} one can show that   $P_j$  are series  convergent  in a neighborhood of $\delta_0$ (see also e.g. \cite{RosR,RSb}).
 Further, by (1.14), we can take
$v_{1+p_n},v_{3+p_n},\cdots,v_{2k-1+p_n}$ as free parameters,
varying near zero. Precisely, if we change them such that
$$0<|v_{1+p_n}|\ll |v_{3+p_n}|\ll\cdots \ll |v_{2k-1+p_n}|\ll 1,
 \ v_{1+p_n}v_{3+p_n}<0,\cdots,v_{2k-1+p_n}v_{2k+1+p_n}<0,$$
 then by (3.7) the function $\bar d$ has exactly  $k$ positive zeros
 in $x_0$ near $x_0=0$, which give $k$ limit cycles. This finishes
 the proof.

By Theorem 1.4 and (3.7) we immediately have

 {\bf Corollary 3.1.} \ {\it Let $(1.9)$ satisfy $(1.11)$ and $(1.12)$
for a  fixed $\delta \in D$. Then, if
 $$v_{2k+1+p_n}(\delta)<0(>0), \ v_{2j-1+p_n}(\delta)=0 \ {\rm for}
 \ j=1,\cdots, k$$
  the origin is a stable $($unstable$)$ focus of order $k$ of $(1.9)$. If
 $$ v_{2j-1+p_n}(\delta)=0 \ {\rm for \ all}
 \ j\geq 1$$the origin is a center of $(1.9).$ }

 {\bf Proof of Theorem
1.7.} \ Under (1.14) $v_{1+p_n},v_{3+p_n},\cdots,v_{2k-1+p_n}$ can
be taken as free parameters.  Further, by our assumption, the origin
is a center of $(1.9)$ as $v_{2j-1+p_n}(\delta)= 0, j=1,\cdots, k.$
It then follows
 $$v_{2j-1+p_n}(\delta)=O(|v_{1+p_n},v_{3+p_n},\cdots,v_{2k-1+p_n}|)
 \ {\rm for \ all } \ j\geq k+1.$$ Therefore, (3.7) can be further
 written in the form
$$\bar d(x_0,\delta)=\sum_{j=1}^kv_{2j-1+p_n}(\delta)x_0^{2j-1+p_n}(1+\bar P_j(x_0,\delta)),$$
where $\bar P_j(0,\delta)=0$  and  $P_j$ are series  convergent  in a neighborhood of $\delta_0$
(\cite{CJ}). Using the reasoning of  Bautin \cite{Bau}
 (see also e.g. \cite{H0,RSb,RosR})
one can easily see that
 the conclusion of the theorem holds. The proof is completed.

{\bf Proof of Theorem 1.8.} \ Now we consider (1.15), where $g$
satisfies (1.11). Let
$$F(x,\delta)=\int_0^xf(x,\delta)dx,\ \
 G(x,\delta)=\int_0^xg(x,\delta)dx.$$ If $f$ satisfies
(1.12), then the origin is a center or focus of (1.15), and
 $$F(\alpha(x,\delta),\delta)-F(x,\delta)=\sum_{j\geq n}B_j(\delta)x^j=\sum_{j\geq n_1}B_j(\delta)x^j,\eqno(3.8)$$
where $$ B_n=\frac{(-1)^n-1}{n}b_{n-1}, \ \ n_1=2l+1, \
l=\left[\frac{n}{2}\right].$$ and $\alpha(x,\delta)=-x+O(x^2)$
 satisfies $G(\alpha(x,\delta),\delta)=G(x,\delta)$ for $|x|$ small.
Note that (1.15) is equivalent to the following system
$$\dot x=y-F(x,\delta), \  \dot y=-g(x,\delta) \eqno(3.9)$$
which has the same Poincar\'e return map $P(x_0,\delta)$ as (1.15).
Introducing the  change of variables $x$ and $t$
$$u=[2nG(x,\delta)]^{\frac{1}{2n}}({\rm sgn }x)=(a_{2n-1})^{\frac{1}{2n}}(x+O(x^2))\equiv \varphi(x),\
\frac{dt}{dt_1}=\frac{\varphi^{2n-1}(x)}{g(x,\delta)}$$ the system
(3.9) becomes
$$\dot u=y-\bar F(u,\delta), \  \dot y=-u^{2n-1},  \eqno(3.10)$$
which is equivalent to
$$\dot u=y, \  \dot y=-u^{2n-1} -y\bar f(u,\delta),\eqno(3.11)$$
where $$\bar F(u,\delta)=F(\varphi^{-1}(u),\delta),\ \bar
f(u,\delta)=\frac{\partial \bar F}{\partial u}(u,\delta).$$ The
systems (3.10) and (3.11) have the same Poincar\'e return map,
denoted by  $P_1(u_0,\delta)$. One can see that the maps $P$ and
$P_1$ have the relation $P_1\circ\varphi=\varphi\circ P$. Hence,
$$P(x_0,\delta)-x_0=K(u_0)(P_1(u_0,\delta)-u_0),$$ where $K(u_0)=(a_{2n-1})^{-\frac{1}{2n}}+O(u_0)$ is analytic. By
(1.13) and (3.7), for $u_0>0$ small we have
$$P_1(u_0,\delta)-u_0=\sum_{j\geq 1}v_j(\delta)u_0^j=\sum_{j\geq 1}v_{2j-1+p_n}(\delta)u_0^{2j-1+p_n}(1+P_j(u_0,\delta)).$$
Hence,
$$\begin{array}{ll}
P(x_0,\delta)-x_0&=\sum_{j\geq
1}v_{2j-1+p_n}(\delta)(a_{2n-1})^{-\frac{1}{2n}}u_0^{2j-1+p_n}(1+\tilde
P_j(u_0,\delta))\\
&=\sum_{j\geq
1}v_{2j-1+p_n}(\delta)(a_{2n-1})^{\frac{2j-2+p_n}{2n}}x_0^{2j-1+p_n}(1+
P_j^*(x_0,\delta)) , \end{array} \eqno(3.12)$$ where $\tilde
P_j(u_0,\delta)=O(u_0)$, $P_j^*(x_0,\delta)=O(x_0)$.

 Since $\alpha$ satisfies $G(\alpha(x,\delta),\delta)=G(x,\delta)$ and $x\alpha<0$  for $|x|$
 small, we have $\varphi(\alpha)=-\varphi(x)$ or $\alpha=\varphi^{-1}(-\varphi(x))=\varphi^{-1}(-u)$,
 where $u=\varphi(x)$. Thus, we have
 $$F(\alpha(x,\delta),\delta)-F(x,\delta)=F(\varphi^{-1}(-u),\delta)-F(\varphi^{-1}(u),\delta)
 =\bar F(-u,\delta)- \bar F(u,\delta).\eqno(3.13)$$
Let $$ \bar f(u,\delta)=\sum_{j\geq n-1}\bar b_j(\delta)u^j.$$ Then
$$\bar F(u,\delta)=\sum_{j\geq n}\frac{\bar b_{j-1}(\delta)}
{j}u^j.$$ Thus, by (3.13) we have
$$F(\alpha(x,\delta),\delta)-F(x,\delta)=-2\sum_{j\geq [n/2]}\frac{\bar b_{2j}(\delta)}
{2j+1}u^{2j+1}.$$ Substituting
$u=\varphi(x)=(a_{2n-1})^{\frac{1}{2n}}(x+O(x^2))$ into the equality
above and comparing with (3.8) we  obtain
$$\begin{array}{l}
B_{2l+1}=-\bar K_l\bar b_{2l}, \  \ B_{2l+2}=O(\bar b_{2l}),\\
B_{2l+2j+1}=-\bar K_{l+j}\bar b_{2l+2j}+O(|\bar b_{2l}, \bar b_{2l+2},\cdots,\bar b_{2l+2j-2}|),\\
B_{2l+2j+2}=O(|\bar b_{2l}, \bar b_{2l+2},\cdots,\bar b_{2l+2j}|), \
j\geq 1,\end{array} \eqno(3.14)$$ where $\bar K_l$, $\bar K_{l+1}$,
$\cdots$ are positive constants.

  Then by  Theorem 1.4  for $u_0>0$ small we clearly have
$$P_1(u_0,\delta)=u_0+\sum_{j\geq 1}v_j(\delta)u_0^j=\sum_{j\geq 1}V_j(\delta)u_0^j,$$
where $V_j$ are introduced before Theorem 1.2.
Thus, by Theorem 1.2, we have
$$v_1=-K_l\bar b_{2l}+O(\bar b_{2l}^2),\ v_{2j+1}|_{v_1=\cdots=v_{2j-1}=0}=-K_{l+j}\bar b_{2l+2j},\ j\geq 1$$
for $n=2l+1$ odd, and
$$v_2=-K_l\bar b_{2l}, \ v_{2j+2}|_{v_2=\cdots=v_{2j}=0}=-K_{l+j}\bar b_{2l+2j},\ j\geq 1$$
for $n=2l$ even, where $K_{l+j}$, $j\geq 0$ are positive constants.
Hence,
$$\begin{array}{l}
v_{1+p_n}=-K_l\bar b_{2l}+(1-p_n)O(\bar b_{2l}^2),\\
v_{2j+1+p_n}=-K_{l+j}\bar b_{2l+2j}+\bar \varphi(\bar b_{2l},\bar
b_{2l+2},\cdots,\bar b_{2l+2j-2}), \ j\geq 1,\end{array}
\eqno(3.15)$$ where $\bar\varphi(0,0,\cdots,0)=0$. Note that (1.15)
is analytic in each $\bar b_j$. It follows from Theorem 1.4 that
$\bar\varphi $ is analytic in $(\bar b_{2l},\bar
b_{2l+2},\cdots,\bar b_{2l+2j-2})$, which yields $\bar\varphi=
O(|\bar b_{2l},\bar b_{2l+2},\cdots,\bar b_{2l+2j-2}|)$. Then (3.14)
and (3.15) together give
$$\begin{array}{l}
v_{1+p_n}=\frac{K_l}{\bar K_l}B_{2l+1}+(1-p_n)O( B_{2l+1}^2),\\
v_{2j+1+p_n}=\frac{K_{l+j}}{\bar K_{l+j}}
B_{2l+2j+1}+O(|B_{2l+1},B_{2l+3},\cdots,B_{2l+2j-1}|), \ j\geq
1,\end{array} \eqno(3.16)$$

Then (1.17) follows from (3.12) and (3.16). This finishes the proof.

{\bf Proof of Theorem 1.9.} \  Let $|\delta-\delta_0|$ be small. For
$n=2$ we have $p_n=1$.   Then the first conclusion is direct from
Corollary 3.1 and (1.17)-(1.20). In fact, we have by Theorem 1.8
$$v_{2k}(\delta_0)=K^*_kB_{2k+1}(\delta_0), \ K^*_k>0, \ v_{2j}(\delta_0)=0 \
{\rm for } \ j=1,\cdots,k-1.$$

 For the second conclusion, we first keep $B_1(\delta)=0$, and vary
 $B_3(\delta),\cdots, B_{2k-1}(\delta)$ near zero to obtain exactly $k-1$
 simple limit cycles near the origin. These limit cycles are bifurcated by
 changing the stability of the focus at the origin $k-1$ times. Then we vary
 $B_1$ such that $0<|B_1|\ll |B_3|$, and $B_1B_3<0$. This step
 produces one more limit cycle bifurcated from the origin by changing the stability of the origin
 which is a node now by \cite{H}. The theorem is proved for the case of $n=2$.

For $n>2$ since $g(-x,\delta)= -g(x,\delta)$,
$f(-x,\delta)=f(x,\delta)$ we have
$$
 b_j(\delta)=0 \ {\rm for} \ j=0,\cdots,n-2\ {\rm and}   \   b_{n-1}^2(\delta)-4na_{2n-1}(\delta)<0$$
if
$$B_{2l+1}(\delta)<0(>0), \ B_{2j-1}(\delta)= 0, j=1,\cdots,
l\eqno(3.17)$$ for some $[n/2]\leq l \leq k$. In this case the
origin is a stable (unstable) focus of (3.9) by Theorem 1.8. If
(3.17) holds for some $0\leq l < [n/2]$, then by \cite{H} again the
origin is a stable (unstable) node of (3.9). Then the proof in this
case is just similar to the above.
 This finishes the proof.

We remark that if $g(-x,\delta)= -g(x,\delta)$ then
$\alpha(x,\delta)=-x$.

{\bf Proof of Theorem 1.10.} \ Consider (1.22). Without loss of
generality, we can assume $X_{m+1}=0$ in (1.22). Otherwise, it needs
only to introduce a change of variables $v=y+X_{m+1}(x,y)$. In this
case, we can write (1.22) into the form
$$\dot x=y,  \ \ \dot y=-g(x)-f(x)y+y^2\sum_{j\geq 0}\varphi_j(x)y^j, \eqno(3.18)$$
where
$$ g(x)=g_m(x)+O(|x|^{m+1}), \ f(x)=f_{m-1}(x)+O(|x|^{m}), \
\varphi_j(x)=O(|x|^{m-1-j}).\eqno(3.19)$$ For the sake of
convenience below, we rewrite the functions $g$, $f$ and $\varphi_j$
as follows:
$$\begin{array}{c}
g(x)=x^{2n-1}[g_0(x)+x^ng_1(x)+x^{2n}g_2(x)+\cdots ],\\
f(x)=x^{n-1}[f_0(x)+x^nf_1(x)+x^{2n}f_2(x)+\cdots] ,\\
\varphi_j(x)=\varphi_{j0}(x)+x^{n-1}\varphi_{j1}(x)+x^{2n-1}\varphi_{j2}(x)+\cdots,
\end{array} \eqno(3.20)$$
where $f_j,\ g_j$ and $\varphi_{jl}$, $j\geq 0,l\geq 1,$ are
polynomials in $x$ of  degree at most $n-1$, and $\varphi_{j0}$,
$j\geq 0$, are polynomials in $x$ with degree at most $n-2$.

Now we change (3.18) by using (2.2) to obtain (2.5) satisfying (2.6)
where
$$\dot x=r^n \sin \theta, \ \ \dot y=r^{2n-1}\sum_{j\geq
0}V_j(\theta,r)r^{jn}, $$ and by (3.20)
$$\begin{array}{cl}
V_0(\theta,r)&=-\cos^{2n-1}\theta g_0(r\cos \theta)-\sin^n\theta
f_0(r\cos \theta)+r\sin^2\theta \varphi_{00}(r\cos \theta),\\
V_j(\theta,r)&=-\cos^{2n-1+jn}\theta g_j(r\cos \theta)-\sin^n\theta
\cos^{jn}\theta f_j(r\cos \theta)+r\sin^{2+j}\theta
\varphi_{j0}(r\cos \theta)\\
&\ \ \  \displaystyle
+\sum_{k=0}^{j-1}\sin^{2+k}\theta\cos^{(j-k)n-1}\theta
\varphi_{k,j-k}(r\cos \theta), \ \ j\geq 1. \end{array} \eqno(3.21)
$$ Hence, we obtain from (2.5)
$$\frac{dr}{d\theta}=r\displaystyle\frac{\sum_{j\geq 0}R_j(\theta,r)r^{jn}}{\sum_{j\geq
0}S_j(\theta,r)r^{jn}},
$$ where
$$\begin{array}{c} S_0(\theta,r)=-n\sin^2\theta+\cos\theta\ V_0(\theta,r), \ R_0(\theta,r)=\sin\theta\cos\theta+\sin \theta\
V_0(\theta,r),\\
S_j(\theta,r)=\cos\theta\ V_j(\theta,r), \ \
R_j(\theta,r)=\sin\theta\ V_j(\theta,r), \  \ j\geq 1. \end{array}
\eqno(3.22)$$ By (3.21) and (3.22) we can further expand $S_j$ and
$R_j$ in $r$ to obtain for $j\geq 0$
$$S_j(\theta,r)=\sum_{l=0}^{n-1}\bar S_{l+jn}(\theta)r^l, \ R_j(\theta,r)=\sum_{l=0}^{n-1}\bar
R_{l+jn}(\theta)r^l \eqno(3.23), $$ so that the above differential
equation can be written as
$$\frac{dr}{d\theta}=r\displaystyle\frac{\sum_{j\geq 0}\bar R_j(\theta)r^{j}}{\sum_{j\geq
0}\bar S_j(\theta)r^{j}}.$$ Further, letting
$$\displaystyle\frac{1}{\sum_{j\geq 0}\bar
S_j(\theta)r^{j}}=\sum_{j\geq 0}\tilde S_j(\theta)r^j$$ and
$$\tilde R_j(\theta)=\sum_{k+l=j}\bar R_k(\theta) \tilde
S_l(\theta), \ j\geq 0 \eqno(3.24)$$
 we  obtain
$$\frac{dr}{d\theta}=r\displaystyle\sum_{j\geq 0}\tilde
R_j(\theta)r^{j}. \eqno(3.25)$$ Note that for  any $ j\geq 0$,
$\tilde S_j$ depends  only  on $\bar S_k$ with $0\leq k\leq j.$ Then
by (3.24) one can see that
 $${\rm For \ any} \ j\geq
0, \ \tilde R_j \ {\rm depends \ only \ on} \ \bar R_k \ {\rm and} \
\bar S_k \ {\rm with} \ 0\leq k\leq j.  \eqno(3.26)$$

Let $r(\theta,r_0)$ denote the solution of (3.25) with the initial
value $r_0$. The for $r_0$ small we have
$$r(\theta, r_0)=\sum_{j\geq 1}r_j(\theta)r_0^j$$
where $r_1,\ r_2, \ r_3, \cdots $ satisfy $r_1(0)=1$,
$r_2(0)=r_3(0)=\cdots=0$, and
$$ \begin{array}{cl}
r'_1&=\tilde R_0r_1,\\
r'_2&=\tilde R_0r_2+\tilde R_1r^2_1,\\
r'_3&=\tilde R_0r_3+2\tilde R_1r_1r_2+\tilde R_2r^3_1,\\
& \cdots \end{array}$$ which implies that for any $j\geq 1$, the
function $r_j$ depends only on $\tilde R_k$ with $ 0\leq k\leq j-1$.
Hence, by Lemma 2.3, (1.23) and (3.26) we come to the following
conclusion:
$$ {\rm For \ any} \ j\geq 1, \ v_j(\delta) \ {\rm depends \ only \ on} \ \bar R_k \ {\rm and} \
\bar S_k \ {\rm with} \ 0\leq k\leq j-1.  \eqno(3.27)$$

Further, by (3.21)--(3.23), one can observe  that for $0\leq l\leq
n-1$, $\bar S_l$ and $\bar R_l$ depend only on the coefficients of
degree $l$ of the polynomials $g_0,\ f_0$ and $x\varphi_{00}$ in
$x$. Hence, by (3.27) we see that for $1\leq j\leq n$, $v_j$ depends
only on the coefficients of degree at most $j-1$ of the polynomials
$g_0,\ f_0$ and $x\varphi_{00}$ in $x$.

Similarly, for $j\geq 1$ and $0\leq l\leq n-1$ or $jn\leq l+jn\leq
(j+1)n-1$, $\bar S_{l+jn}$ and $\bar R_{l+jn}$ depend only on the
coefficients of degree $l$ of the polynomials $g_j,\ f_j$,
$x\varphi_{j0}$ and $\varphi_{i,j-i}$ with $i=0,\cdots,j-1$ in $x$.
In other words,  for  $jn+1\leq u\leq (j+1)n$, $\bar S_{u-1}$ and
$\bar R_{u-1}$ depend only on the coefficients of degree $u-1-jn$ of
the polynomials $g_j,\ f_j$, $x\varphi_{j0}$ and $\varphi_{i,j-i}$
with $i=0,\cdots,j-1$ in $x$. Let $N_{[a,b]}$ denote the set of
integers in the interval $[a,b]$. Then for $jn+1\leq u\leq (j+1)n$,
we have
$$N_{[0,u-1]}=\bigcup_{i=0}^{j-1}N_{[in,(i+1)n-1]}\bigcup
N_{[jn,u-1]}.$$ Thus, for all $k\in N_{[in,(i+1)n-1]}$, $\bar S_{k}$
and $\bar R_{k}$ depend only on $g_i,\ f_i$, $x\varphi_{i0}$ and
$\varphi_{l,i-l}$ with $l=0,\cdots,i-1$. And for $k\in
N_{[jn,u-1]}$, $\bar S_{k}$ and $\bar R_{k}$ depend only on the
coefficients of degree $k-jn$ of the polynomials $g_j,\ f_j$,
$x\varphi_{j0}$ and $\varphi_{l,j-l}$ with $l=0,\cdots,j-1$ in $x$.

Therefore, by (3.27) for $jn+1\leq u\leq (j+1)n$, $v_u(\delta)$
depends only on the functions $g_i,\ f_i$, $x\varphi_{i0}$ and
$\varphi_{l,i-l}$ with $l=0,\cdots,i-1$, $i=0,\cdots,j-1$ and the
coefficients of degree at most $u-1-jn$ of the polynomials $g_j,\
f_j$, $x\varphi_{j0}$ and $\varphi_{l,j-l}$ with $l=0,\cdots,j-1$ in
$x$.

We claim that if $j\geq 0$, $m\geq (j+1)n$, then for $jn+1\leq u\leq
(j+1)n$, $v_u(\delta)$ depends only on the functions $g_i,\ f_i$,
 with $i=0,\cdots,j-1$ and the coefficients of degree at most $u-1-jn$ of
the polynomials $g_j,\ f_j$
 in $x$.

In fact, by the above discussion, we need only to prove
$\varphi_{00}=0$ in the case $j=0$ and $\varphi_{ls}=0$ for $l+s\leq
j$ and $0\leq l\leq j-1$ in the case $j>0$. This can be shown easily
since $$\varphi_{j0}=O(|x|^{m-1-j}), \ \varphi_{js}=O(|x|^{m-j-sn})
\ {\rm for} \ s\geq 1$$ and $$\deg \varphi_{j0}\leq n-2, \ \ \deg
\varphi_{js} \leq n-1 \ {\rm for} \ s\geq 1$$ by (3.19) and (3.20).

By (3.20) again, the above claim can be restated that if $j\geq 0$,
$m\geq (j+1)n$, then for $jn+1\leq u\leq (j+1)n$, $v_u(\delta)$
depends only on the coefficients of degree at most $2n+u-2$ of $g$
and   the coefficients of degree at most $n+u-2$ of $f$
 in $x$.
Thus, for any integers $k$ and $m$ satisfying  $k\geq 1$ and $m\geq
(k+1)n$, by taking $j=0,\cdots,k$ we know that for all $1\leq
u\leq(k+1)n$, $v_u(\delta)$ depends only on the coefficients of
degree at most $2n+u-2$ of $g$ and   the coefficients of degree at
most $n+u-2$ of $f$
 in $x$.

Finally, by (3.19), if $m\geq (k+3)n-2$ then
$$2n+u-2\leq m, \ \ n+u-2\leq m-1 \ \ {\rm for} \ u\leq (k+1)n. $$In
this case, for all $1\leq u\leq(k+1)n$, $v_u(\delta)$ depends only
on $g_m$ and $f_{m-1}$ in (3.19). Then the conclusion of Theorem
1.10 follows.

\section{Application examples}

In this section  we give some application examples based on the
examples given in \cite{AG2}.
 Consider a Kukles type system of the form
 $$\dot x=y, \ \dot y=-( a_{11}xy+a_{02}y^2+a_{30}x^3+
 a_{21}x^2y+a_{12}xy^2+a_{03}y^3).\eqno(4.1)$$
The authors \cite{AG2} proved that if $a_{30}>0$ and
$a^2_{11}-8a_{30}<0$ then for (4.1) $v_2=v_4=v_6=v_8=0$ if and only
if $a_{21}=a_{03}=a_{11}a_{02}=0$, which implies that the origin is
a center. Moreover, there can be  3 limit cycles near the origin.
See Theorem 4.1 in \cite{AG2} and its proof.

Based on this conclusion and by Theorem 1.6 we have immediately

{\bf Proposition 4.1.} \ {\it Let $a_{11},a_{02},a_{30},
 a_{21},a_{12}$ and $a_{03}$ be bounded parameters satisfying
 $$ a_{30}>0, \ a^2_{11}-8a_{30}<0, \ |a_{21}|+|a_{03}|+|a_{11}a_{02}|>0.$$
  Then there exists a neighborhood $V$ of the
 origin such that
 the system $(4.1)$ has at most $3$ limit cycles in
 $V$.}

Then consider
$$\dot x=-y+Ax^2+Bxy+Cy^2,\ \dot y=x^3+xy^2+y^3.\eqno(4.2)$$
By Theorem 4.2 in \cite{AG2} and its proof if $A^2<2$ then the
origin of (4.2) is always a focus with $|v_2|+|v_4|+|v_6|+|v_8|>0$.
Moreover, there are systems inside (4.2) with at least 3 limit
cycles around the origin. Then by Theorem 1.6 again we have

{\bf Proposition 4.2.} \ {\it Let $A,B$ and $C$ be bounded
parameters with $A^2<2$. Then   there exists a neighborhood $V$ of
the origin such that  the system $(4.2)$ has at most $3$ limit
cycles in
 $V$.}

Finally, consider
$$\dot x=y,\ \dot y=-(x^3+x^5)-\sum_{j=0}^kb_{2j}x^{2j}y,\eqno(4.3)$$
where $k\geq 2.$ By Theorems 1.7-1.9, we obtain

{\bf Proposition 4.3.} \ {\it Let $b_{2j}$ be bounded parameters.
Then

$(1)$ \ If $b_0=0$,  the system $(4.3)$ has at most $k-1$ limit
cycles near the origin; and $k-1$ limit cycles can appear.

$(2)$ \ If $b_0\neq 0$,  there are systems inside $(4.3)$ which have
at least $k$ limit cycles near the origin. }

%\section*{Acknowledgements}

\end{document}